\documentclass{amsart}

\usepackage{graphicx}
\usepackage[pdfauthor="Richard J. Mathar",colorlinks=true,bookmarks=true,bookmarksopen=true,citecolor=blue]{hyperref}
\newtheorem{defn}{Definition}
\newtheorem{rem}{Remark}
\newtheorem{exa}{Example}

\begin{document}

\title[Chebyshev expansion of $x^m(-\log x)^l$]{Chebyshev approximation of $x^m(-\log x)^l$ in the interval $0\le x\le 1$}

\author{Richard J. Mathar}
\urladdr{https://www.mpia-hd.mpg.de/homes/mathar}
\email{mathar@mpia.de}
\address{Max-Planck Institute for Astronomy, K\"onigstuhl 17, 69117 Heidelberg, Germany}


\date{\today}

\begin{abstract}
The series expansion of $x^m(-\log x)^l$ in terms of the shifted Chebyshev Polynomials $T_n^*(x)$ requires
evaluation of the integral family $\int_0^1 x^m (-\log x)^l dx/\sqrt(x-x^2)$. We demonstrate that these
can be reduced by partial integration to sums over integrals with exponent $m=0$ which have known 
representations as finite sums over polygamma functions.
\end{abstract}

\maketitle

\section{Chebyshev Series Coefficients}
The shifted Chebyshev Polynomials of the first kind are defined as \cite{Clenshaw1954}\cite[22.5.5]{AS}
\begin{defn} (Shifted Chebyshev Polynomials)
\begin{equation}
T^*_n(x) = \cos(n\theta) ;\quad \cos\theta =2x-1.
\end{equation}
\end{defn}
All properties follow from their representation
as Chebyshev Polynomials $T_n(x)$ of the first kind,
\begin{equation}
T_n(u) =  T_n^*((1+u)/2); \quad
T_n(2x-1) =  T_n^*(x). \quad
\end{equation}
They are terminating Gaussian Hypergeometric Functions \cite[8.942.1]{GR}:
\begin{equation}
T_n^*(x) = {}_2F_1(n,-n;1/2;1-x).
\end{equation}
Efficient evaluation of sums over these polynomials is generally using the
recurrence \cite[Tab. 18.9.1]{DLMF}
\begin{equation}
T^*_{n+1}(x) = [(4-2\delta_{n,0})x-2+\delta_{n,0}] T^*_n(x)-T^*_{n-1}(x).
\label{eq.Tstarrec}
\end{equation}

\begin{exa}\label{exa.Tstar}
\begin{eqnarray}
T^*_0(x)&=& 1 ; \\
T^*_1(x)&=& 2x-1 ; \label{eq.Tstar1}\\
T^*_2(x)&=& 8x^2-8x+1 ;\\
T^*_3(x)&=& 32x^3-48x^2+18x-1.
\end{eqnarray}
See the table \cite[A127674]{sloane} in the Online Encyclopedia of Integer Sequences (OEIS) for more of these integer coefficients.
\end{exa}

The orthogonality relation is \cite[18.3.1]{DLMF}
\begin{equation}
\int_0^1 \frac{T_n^*(x)T_m^*(x)}{\sqrt{x-x^2}}dx = \left\{
\begin{array}{ll}
\delta_{n,m}\pi,& n=0;\\
\delta_{n,m}\frac{\pi}{2},& n>0.\\
\end{array}
\right.
\label{eq.ortho}
\end{equation}

\begin{defn} (Expansion coefficients)
The expansion coefficients $A_{.,.,.}$ are defined as
\begin{equation}
x^m(-\log x)^l\equiv \mathop{{\sum}'}_{n= 0}^\infty A_{m,l,n}T_n^*(x)
\label{eq.Adef}
\end{equation}
where the prime at the sum symbol specifies that the term at $n=0$ contributes
halved to the sum.
\end{defn}
Projection of the expansion coefficients with the aid of the orthogonality 
describes the $A$ as elementary integrals
\begin{equation}
A_{m,l,k} 
=
\frac{2}{\pi}
\int_0^1 x^m(-\log x)^l \frac{T_k^*(x)}{\sqrt{x-x^2}} dx, \quad k=0,1,2,\ldots
\end{equation}
This manuscript is a guide to the numerical evaluation of these integrals.
Convergence properties are \emph{not} discussed \cite{BoydAMC29}.
If the $T_k^*$-factor is absent \cite[4.272.16]{GR}\cite{QiRACSAM114},
\begin{multline}
A_{m,l,0} 
=
\frac{2}{\pi}\int_0^1 \frac{x^{m-1/2}(-\log x)^l}{\sqrt{1-x}}dx
=
\frac{2}{\pi}l!\sum_{k\ge 0}\binom{-1/2}{k}\frac{(-)^k}{(m+1/2+k)^{l+1}}
\\
=
\frac{2}{\pi}l!\sum_{k\ge 0}\frac{(2k-1)!!}{(2k)!!}\frac{1}{(m+1/2+k)^{l+1}}
=
\frac{2}{\pi}l!\sum_{k\ge 0}\binom{2k}{k}\frac{1}{4^k(m+1/2+k)^{l+1}}
\\
=
\frac{2}{\pi}\frac{l!}{(m+1/2)^{l+1}}
{}_{l+2}F_{l+1}\left(\begin{array}{cccc}
1/2 & m+1/2 & m+1/2... & m+1/2 \\
& m+3/2 & m+3/2... & m+3/2 \\
\end{array}\mid 1\right)
.
\label{eq.AasF}
\end{multline}
\begin{defn} (Double Factorial)
\begin{equation}
u!! = \left\{ 
\begin{array}{ll}
u(u-2)(u-4)\cdots 4\cdot 2; & u\, \mathrm{even}; \\
u(u-2)(u-4)\cdots 3\cdot 1; & u\, \mathrm{odd}.
\end{array}
\right.
\end{equation}
\end{defn}
\begin{rem}
Numerical values of these $_{l+2}F_{l+1}$-constants at $m=0$, $l=0$--$3$,
are sequences 
A019669, A173623, A318742, and A375594
in the OEIS \cite{sloane}.
\end{rem}

\section{No Logarithm, $l=0$}
A somewhat degenerate case is a kernel without the logarithm at $l=0$; it basically expands $x^m$ in a series of 
shifted Chebyshev polynomials, and inverts the triangular matrix of coefficients implied in Example \ref{exa.Tstar}
that represents $T^*$ as an integer matrix times the monomial basis:

\begin{equation}
\left(
\begin{array}{c}
T_0^* \\
T_1^* \\
T_2^* \\
T_3^* \\
T_4^* \\
\ldots
\end{array}
\right)
=
\left(
\begin{array}{ccccc}
1 & 0 & 0 & 0 & \ldots \\
-1 & 2 & 0 & 0 & \ldots \\
1 & -8 & 8 & 0 & \ldots \\
-1 & 18 & -48 & 32 & \ldots \\
1 & 32 & 160 & -256 & 128 \\
\ldots
\end{array}
\right)
\cdot
\left(
\begin{array}{c}
1 \\
x \\
x^2 \\
x^3 \\
x^4 \\
\ldots
\end{array}
\right)
\end{equation}
\begin{equation}
\leftrightarrow 
\left(
\begin{array}{c}
1 \\
x \\
x^2 \\
x^3 \\
\ldots
\end{array}
\right)
=
\left(
\begin{array}{ccccc}
1 & 0 & 0 & 0 & \ldots \\
1/2 & 1/2 & 0 & 0 & \ldots \\
3/8 & 1/2 & 1/8 & 0 & \ldots \\
5/16 & 15/32 & 3/16 & 1/32 & \ldots \\
35/128 & 7/16 & 7/32 & 1/16 & 1/128 
\ldots
\end{array}
\right)
\cdot
\left(
\begin{array}{c}
T_0^* \\
T_1^* \\
T_2^* \\
T_3^* \\
\ldots
\end{array}
\right).
\label{eq.Tinv}
\end{equation}

An immediate application of \cite[2.261]{GR} is
\begin{equation}
A_{0,0,0} 
= 
\frac{2}{\pi}
\int_0^1 \frac{1}{\sqrt{x-x^2}}dx
= 
\frac{2}{\pi}\pi
=2
.
\end{equation}
Higher $m$ are then recursively covered by \cite[2.263.1]{GR}
\begin{equation}
A_{m,0,0} 
= 
\frac{2}{\pi}
\int_0^1 x^m \frac{1}{\sqrt{x-x^2}}dx
= 
\frac{2}{\pi}
[
\frac{(2m-1)}{2m}\int_0^1 \frac{x^{m-1}}{\sqrt{x-x^2}}dx 
]
= 
\frac{2m-1}{2m} A_{m-1,0,0}.
\end{equation}

In summary
\begin{equation}
A_{m,0,0} 
=2\frac{(2m-1)!!}{(2m)!!}
=\frac{(2m-1)!!}{2^{m-1}m!}
\label{eq.l0}
\end{equation}
are twice the entries in the left column of the matrix in \eqref{eq.Tinv}---the prime at
the sum symbol in \eqref{eq.Adef} established the division by 2.

The representation \eqref{eq.Tstar1} leads to
\begin{equation}
A_{m,l,1}  = 2 A_{m+1,l,0}-A_{m,l,0}.
\label{eq.recn1}
\end{equation}
Likewise \eqref{eq.Tstarrec} yields
\begin{equation}
A_{m,l,n}  = 4 A_{m+1,l,n-1}-2A_{m,l,n-1}-A_{m,l,n-2}, \quad n\ge 2.
\label{eq.recn2}
\end{equation}

\section{Logarithm to the first power, $l=1$}
Two special functions in this business are \cite[8.384,8.36]{GR}
\begin{defn} (Beta function)
\begin{equation}
B(x,y) = \frac{\Gamma(x)\Gamma(y)}{\Gamma(x+y)}.
\end{equation}
\end{defn}
\begin{defn} (Digamma function)
\begin{equation}
\psi(x) = \frac{\Gamma'(x)}{\Gamma(x)}
\end{equation}
\end{defn}
The substitution $x=1-u$ generates
\begin{equation}
A_{m,1,0} 
=
\frac{2}{\pi}
\int_0^1 x^m(-\log x) \frac{1}{\sqrt{x-x^2}} dx
=
 \frac{2}{\pi}
(-)\int_0^1 (1-u)^m \log (1-u) \frac{1}{\sqrt{u-u^2}} du
,
\end{equation}
which is found in the Gradsteyn-Ryzhik tables \cite[4.293.13,4.256]{GR}:
\begin{multline}
A_{m,1,0} 
=
 \frac{2}{\pi}
(-)\int_0^1 (1-u)^{m+1/2-1}u^{1/2-1} \log (1-u) du
\\
=
 \frac{2}{\pi}
(-)B\left(\frac12,m+\frac12\right)\left[\psi(m+1/2)-\psi(m+1)\right]
\\
=
- 
\frac{(2m-1)!!}{2^{m-1}m!}[\psi(m+1/2)-\psi(m+1)]
.
\end{multline}
\begin{rem}
The same result is obtained by applying
\cite{GottschalkJPA21}
\begin{equation}
_{p+1}F_p\left(\begin{array}{c} a,a,\ldots,b \\ a+1,a+1,\ldots, a+1
\end{array}
\mid
 1\right)
=
-\frac{\Gamma(1-b)(-a)^p}{\Gamma(p)}\left[\frac{d^{(p-1)}}{dr^{(p-1)}}\frac{\Gamma(r)}{\Gamma(r+1-b)}\right]_{r=a}, b\notin \mathbb{Z}^+
\end{equation}
to \eqref{eq.AasF}.
\end{rem}
By well-known properties of the digamma function \cite[6.3.8,6.3.6]{AS}\cite{BeumerAMM68}\cite[8.370,8.375.2]{GR}
\begin{equation}
A_{m,1,0}=
- 
\frac{(2m-1)!!}{2^{m-1}m!}\left[\sum_{i=0}^{m-1}\frac{1}{(i+1)(2i+1)}-2\log 2\right]
.
\label{eq.l1}
\end{equation}

\begin{exa}
\begin{eqnarray}
A_{0,1,0} &=& 4\log 2 \approx 2.77258872 \label{eq.A010};\\
A_{1,1,0} &=& -1+2\log 2 \approx 0.386294 ;\\
A_{2,1,0} &=& -\frac{7}{8}+\frac{3}{2}\log 2 \approx 0.16472077 ;\\
A_{3,1,0} &=& -\frac{37}{48}+\frac{5}{4}\log 2 \approx 0.095600 ;\\
A_{4,1,0} &=& -\frac{533}{768}+\frac{35}{32}\log 2 \approx 0.064119 .
\end{eqnarray}
\end{exa}
\begin{exa}
The property \eqref{eq.recn1} leads to
\begin{eqnarray}
A_{1,1,1} &=& -\frac34 +\log 2 \approx -0.056852819440 ;\\
A_{2,1,1} &=& -\frac23 +\log 2 \approx 0.0264805138932 ;\\
A_{3,1,1} &=& -\frac{79}{128} +\frac{15}{16}\log 2 \approx 0.03263798177 ;\\
A_{4,1,1} &=& -\frac{277}{480} +\frac{7}{8}\log 2 \approx 0.02942044965 .
\end{eqnarray}
\end{exa}
\begin{exa}
Application of \eqref{eq.recn2} yields for example
\begin{eqnarray}
A_{1,1,3} &=& \frac{1}{24} \approx 0.04166666666\ldots ;\\
A_{2,1,3} &=& -\frac{1}{40} = -0.025 ;\\
A_{3,1,3} &=& -\frac{49}{640}+\frac{1}{16}\log 2 \approx -0.0332408012150034\ldots;\\
A_{4,1,3} &=& -\frac{129}{1120}+\frac{1}{8}\log 2 \approx -0.028535173858578 .
\end{eqnarray}
\end{exa}

Decrementing $n$ down to zero with the aid of \eqref{eq.recn1}, \eqref{eq.recn2} and \eqref{eq.l1} 
reduces all $A_{m,1,n}$ to sums of rational numbers and rational multiples of $\log 2$.

\section{Logarithmic powers $l\ge 2$} 
\subsection{m=0}\label{sec.m0}

The substitution $x=t^2$
shows that 
---apart from a factor $2/\pi$--the $A_{0,l,0}$ are the log-sine integrals
\begin{equation}
\int_0^1 \frac{(-\log x)^l}{\sqrt{x-x^2}} dx
=
2 (-)^l \int_0^1 \frac{(\log t^2)^l}{\sqrt{1-t^2}} dt
=
2 (-2)^l \int_0^1 \frac{(\log t)^l}{\sqrt{1-t^2}} dt
=
2 (-2)^l \int_0^{\pi /2} (\log \sin \theta )^l d\theta
\label{eq.xt2}
\end{equation}
studied earlier
\cite{KolbigMathComp40,BowmanJLMS22,ChaudhuriAMM74}:

\begin{gather}
\int_0^1 \frac{ \log (x)}{\sqrt{x-x^2}}dx= -2\pi \log 2 \approx -4.355172...; 
\\
\int_0^1 \frac{ \log^2 (x)}{\sqrt{x-x^2}}dx= \frac{\pi^3}{3}+4\pi \log^2 2\approx 16.3729762...;
\\
\int_0^1 \frac{ \log^3 (x)}{\sqrt{x-x^2}}dx= -12\pi\zeta(3)-2\pi^3 \log 2 -8\pi \log ^3 (2) \approx -96.6701265...;
\\
\int_0^1 \frac{ \log^4 (x)}{\sqrt{x-x^2}}dx= 8\pi^3\log^2 2+96\pi\log 2 \zeta(3)+16\pi \log^4 2 +\frac{19}{15}\pi^5 \approx 769.692354...;
\\
\int_0^1 \frac{ \log^5 (x)}{\sqrt{x-x^2}}dx= 
-480\pi \zeta(3)\log^2 2 -32\pi \log^5 2-720\pi \zeta(5) -\frac{38}{3}\pi^5 \log 2 
\\
-\frac{80}{3} \pi^3\log^32 -40 \pi^3\zeta(3)\approx -7685.47786; \nonumber
\\
\int_0^1 \frac{ \log^6 (x)}{\sqrt{x-x^2}}dx= 
76\pi^5 \log^2 2+1440\pi \zeta^2(3)+\frac{275}{21}\pi^7+64\pi \log^62 \\
+480 \zeta(3)\pi^3\log 2 
+1920 \pi \zeta(3)\log^3 2
+8640\pi \zeta(5)\log 2 +80 \pi^3 \log^4 2 \approx 92181.5543; \nonumber
\\
\int_0^1 \frac{ \log^7 (x)}{\sqrt{x-x^2}}dx= 
-3360\pi^3\zeta(3)\log^22 -532\pi^5 \zeta(3)-5040\pi^3\zeta(5)-128\pi \log^7 2 
\\
-\frac{1064}{3}\pi^5\log^3 2
-20160\pi \zeta^2(3)\log 2
-6720\pi \zeta(3)\log^4 2 -224\pi^3\log^5 2 -\frac{550}{3}\pi^7\log 2
\nonumber
\\
-60480 \pi \zeta(5) \log^22 -90720\pi \zeta(7) 
\approx -0.12903396\times 10^7 ;
\nonumber
\\
\int_0^1 \frac{ \log^8 (x)}{\sqrt{x-x^2}}dx= 
17920\pi^3\zeta(3)\log^3 2
+\frac{4400}{3}\pi^7\log^2 2
+\frac{1792}{3}\pi^3\log^6 2
+\frac{4256}{3}\pi^5\log^4 2
\\
+1451520\pi\zeta(7)\log 2
+21504\pi\zeta(3)\log^5 2
+161280\pi\zeta^2(3)\log^2 2
\nonumber
\\
+322560\pi\zeta(5)\log^3 2
+256\pi\log^8 2
+80640\pi^3\zeta(5)\log 2
+\frac{11813}{45}\pi^9
\nonumber
\\
+8512\pi^5\zeta(3)\log 2
+13440\pi^3\zeta^2(3)
+483840\pi\zeta(3)\zeta(5) \approx 0.20644368..\times 10^8.
\nonumber
\end{gather}

For larger $l$ the series representation
\eqref{eq.AasF}
converge quickly to establish the values for numerical purposes.
The hypergeometric series \eqref{eq.AasF} leads to representations that are essentially Bell-polynomials
of polygamma-functions \cite{KrupnikovJCAM78}.
\begin{rem}
Adamchik writes the integral as
\cite[(16)]{AdamchikJCAM79}
\begin{equation}
A_{0,l,0} 
= 
(-)^l \frac{2}{\pi}\int_0^1 x^{1/2-1}\frac{\log ^lx }{(1-x)^{1/2}} dx
= 
(-)^l \frac{2}{\pi}\Gamma(1/2)\Gamma(l+1)\left[\begin{array}{c}1/2\\ l+1\end{array}\right]
\label{eq.A12}
\end{equation}
in terms of generalized Stirling numbers of the First Kind
or generalized hypergeometric series
\begin{equation}
\left[\begin{array}{c}z\\ p\end{array}\right]
\equiv \frac{(-1)^{p+1}}{z^p\Gamma(1-z)}{}_{p+1}F_p\left(\begin{array}{c}z, z, \ldots , z
\\ z+1, \ldots, z+1\end{array}\mid 1\right)
=\frac{(-)^{p+1}\sin(\pi z)}{\pi}\sum_{k\ge 0} \frac{\Gamma(k+z)}{k!(k+z)^p}
.
\end{equation}
\end{rem}

\subsection{m$>$0}
We continue for $m>0$ with the format \eqref{eq.xt2}
\begin{equation}
A_{m,l,0} 
= 
\frac{4}{\pi}(-2)^l
\int _0^1 t^{2m} \frac{(\log t)^l}{\sqrt{1-t^2}} dt
\end{equation}
and its partial integration
\begin{multline}
= 
\frac{4}{\pi}(-2)^l
x^{2m} \frac{\log^{l-1} x}{\sqrt{1-x^2}}x[\log x-1]\mid_0^1
\\
-
\frac{4}{\pi}(-2)^l
\int_0^1 dx 
x[\log x-1]
\left[2mx^{2m-1}\frac{\log^{l-1} x}{\sqrt{1-x^2}}
+(l-1)x^{2m-1}\frac{\log^{l-2}x}{\sqrt{1-x^2}}
+x^{2m+1}\frac{\log^{l-1} x}{(1-x^2)^{3/2}}
\right]
.
\end{multline}
Assuming $2m\ge 1$ and $l\ge 2$ the first term on the right hand side vanishes:
\begin{multline}
A_{m,l,0} 
 =
-
\frac{4}{\pi}(-2)^l
\int_0^1 dx
\left[2mx^{2m}\frac{\log^{l} x}{\sqrt{1-x^2}}
+(l-1)x^{2m}\frac{\log^{l-1}x}{\sqrt{1-x^2}}
+x^{2m+2}\frac{\log^{l} x}{(1-x^2)^{3/2}}
\right]
\\
+
\frac{4}{\pi}(-2)^l
\int_0^1 dx
\left[2mx^{2m}\frac{\log^{l-1} x}{\sqrt{1-x^2}}
+(l-1)x^{2m}\frac{\log^{l-2}x}{\sqrt{1-x^2}}
+x^{2m+2}\frac{\log^{l-1} x}{(1-x^2)^{3/2}}
\right]
\\
=
-2mA_{m,l,0}
-
\frac{4}{\pi}(-2)^l
(l-1)
\int_0^1 dx
x^{2m}\frac{\log^{l-1}x}{\sqrt{1-x^2}}
-
\frac{4}{\pi}(-2)^l
\int_0^1 dx
x^{2m+2}\frac{\log^{l} x}{(1-x^2)^{3/2}}
\\
+
\frac{4}{\pi}(-2)^l
2m
\int_0^1 dx
x^{2m}\frac{\log^{l-1} x}{\sqrt{1-x^2}}
+
\frac{4}{\pi}(-2)^l
(l-1)
\int_0^1 dx
x^{2m}\frac{\log^{l-2}x}{\sqrt{1-x^2}}
+
\frac{4}{\pi}(-2)^l
\int_0^1 dx
x^{2m+2}\frac{\log^{l-1} x}{(1-x^2)^{3/2}}
\\
=
-2mA_{m,l,0}
+
2
\frac{4}{\pi}(-2)^{l-1}
(l-1-2m)
\int_0^1 dx
x^{2m}\frac{\log^{l-1}x}{\sqrt{1-x^2}}
-
\frac{4}{\pi}(-2)^l
\int_0^1 dx
x^{2m+2}\frac{\log^{l} x}{(1-x^2)^{3/2}}
\\
+
4
\frac{4}{\pi}(-2)^{l-2}
(l-1)
\int_0^1 dx
x^{2m}\frac{\log^{l-2}x}{\sqrt{1-x^2}}
-2
\frac{4}{\pi}(-2)^{l-1}
\int_0^1 dx
x^{2m+2}\frac{\log^{l-1} x}{(1-x^2)^{3/2}}
.
\label{eq.Pint}
\end{multline}
To keep an overall view of that recurrence, define a family of integrals
with 
a $3/2$-power
in the denominator:
\begin{defn} (Companion Integrals)
\begin{equation}
B_{m,l}\equiv \frac{4}{\pi}(-2)^l \int_0^1 x^{2m+2}\frac{\log^l (x)}{(1-x^2)^{3/2}} dx
.
\label{eq.Bdef}
\end{equation}
\end{defn}
A series representation is \cite[4.272.16]{GR}
\begin{multline}
B_{m,l}
= \frac{2}{\pi}\int_0^1 x^{m+1/2} \frac{(-\log x) ^l}{(1-x)^{3/2}} dx
\\
=
\frac{2}{\pi}l!\sum_{k\ge 0}\binom{-3/2}{k} \frac{(-)^k}{(m+3/2+k)^{l+1}}
=
\frac{2}{\pi}l!\sum_{k\ge 0}\frac{(2k+1)!!}{(2k)!!}\frac{1}{(m+3/2+k)^{l+1}}
\\
=
\frac{2l!}{\pi(m+3/2)^{l+1}} {}_{l+2}F_{l+1}\left(\begin{array}{cccc} 
3/2, m+3/2, m+3/2, ..., m+3/2\\
m+5/2, m+5/2, ..., m+5/2
\end{array}\mid 1\right)
.
\label{eq.BasF}
\end{multline}
In \eqref{eq.Bdef} the denominator produces a non-lifted pole at $x=1$ if $l=0$:
\begin{equation}
B_{m,0}=\frac{4}{\pi} \int_0^1 x^{2m+2}/(1-x^2)^{3/2} dx = \infty ;
\end{equation}

\begin{exa}
Eq. \eqref{eq.BasF} turns into a well-known special case 
of the Gaussian Hypergeometric Function of unit argument if $m=-1$ and $l=1$ \cite[15.1.20]{AS}:
\begin{multline}
B_{-1,1}
=
\frac{2}{\pi(1/2)^2} {}_{3}F_{2}\left(\begin{array}{c} 
3/2, 1/2, 1/2\\
3/2, 3/2
\end{array}\mid 1\right)
=
\frac{2}{\pi(1/2)^2} {}_{2}F_{1}\left(\begin{array}{c} 
1/2, 1/2\\
3/2
\end{array}\mid 1\right)
=4.
\end{multline}
\end{exa}
\begin{exa}
Eq. \eqref{eq.BasF} 
allows a reverse lookup via \eqref{eq.AasF} when $m=-1$:
\begin{multline}
B_{-1,l} 
= \frac{4}{\pi}(-2)^l \int_0^1  \frac{\log^l x}{(1-x^2)^{3/2}}dx
=
\frac{2l!}{\pi(1/2)^{l+1}} {}_{l+2}F_{l+1}\left(\begin{array}{c} 
3/2, 1/2, 1/2,\ldots, 1/2\\
3/2, 3/2, \ldots,3/2
\end{array}\mid 1\right)
\\
=
\frac{2^{l+2}l!}{\pi} {}_{l+1}F_{l}\left(\begin{array}{c} 
 1/2, 1/2,\ldots, 1/2\\
3/2, \ldots,3/2
\end{array}\mid 1\right)
=
\frac{2^{l+2}l!}{\pi} \frac{A_{0,l-1,0}}{\frac{2}{\pi} \frac{(l-1)!}{(1/2)^l}}
 =
2lA_{0,l-1,0},\quad l \ge 1.
\label{eq.Bm1}
\end{multline}
\end{exa}
Lowering the first index of $B$ happens with
\begin{multline}
B_{m+1,l}
=\frac{4}{\pi} (-2)^l\int_0^1 x^{2m+4}\frac{\log^l(x)}{(1-x^2)^{3/2}} dx 
\\
=
-\frac{4}{\pi} (-2)^l\int_0^1 x^{2m+2}(1-x^2-1)\frac{\log^l(x)}{(1-x^2)^{3/2}} dx 
=
-A_{m+1,l,0} +B_{m,l}.
\label{eq.Bmrec}
\end{multline}
Terminating at \eqref{eq.Bm1} this is
\begin{equation}
B_{m,l}
=
-\sum_{s=0}^m A_{s,l,0} +2lA_{0,l-1,0}.
\end{equation}
\begin{rem}
In Adamchik's notation $m=0$ refers to a generalized Stirling number \cite{AdamchikJCAM79}
\begin{multline}
B_{0,l} 
= \frac{2}{\pi}(-)^l \int_0^1 x^{1/2} \frac{\log ^l x}{(1-x)^{3/2}} dx.
= \frac{2}{\pi}(-)^l \Gamma(-1/2)\Gamma(l+1)
\left[\begin{array}{c}3/2\\ l+1\end{array}\right].
\end{multline}
The recurrence  \cite[24.1.3]{AS}
\begin{equation}
\left[\begin{array}{c}n\\ k\end{array}\right]
=
(n-1)
\left[\begin{array}{c}n-1\\ k\end{array}\right]
+
\left[\begin{array}{c}n-1\\ k-1\end{array}\right]
\end{equation}
for these Stirling numbers and \eqref{eq.A12} establish
\begin{equation}
B_{0,l} 
= 
-A_{0,l,0}
+ 
2l A_{0,l-1,0}
.
\end{equation}
This follows also by inserting \eqref{eq.Bm1} into the right hand side of \eqref{eq.Bmrec}.
\end{rem}

With these companion integrals \eqref{eq.Pint} may  be written as
\begin{equation}
(1+2m)A_{m,l,0}= 2(l-1-2m)A_{m,l-1,0}-B_{m,l} +4(l-1)A_{m,l-2,0} -2B_{m,l-1}.
\label{eq.tmp1}
\end{equation}

Replacing the two $B$-terms 
with the aid of \eqref{eq.Bmrec} yields
\begin{multline}
2m A_{m,l,0}= 2(l-2m)A_{m,l-1,0}-B_{m-1,l} +4(l-1)A_{m,l-2,0} -2B_{m-1,l-1}
\\
= 
\sum_{s=0}^{m-1} A_{s,l,0} 
+2(l-2m)A_{m,l-1,0}
-2lA_{0,l-1,0}
+2\sum_{s=0}^{m-1} A_{s,l-1,0}
+4(l-1)[A_{m,l-2,0} 
-A_{0,l-2,0}].
\label{eq.recm}
\end{multline}
For $m>0$ this equation 
reduces every  $A_{m,l,0}$ to a finite sum of rational numbers, rational multiples of $\log 2$
and rational multiples of $A_{0,l',0}$, $l'\le l$.
\begin{exa}
\begin{eqnarray}
A_{3,3,0} &=& -\frac{769}{288}-\frac{155}{24}\log 2 -\frac{37}{32}A_{0,2,0}+\frac{5}{16}A_{0,3,0}; \\
A_{4,3,0} &=& -\frac{40673}{18432}-\frac{4163}{768}\log 2 -\frac{533}{512}A_{0,2,0}+\frac{35}{128}A_{0,3,0}; \\
A_{2,4,0} &=& -\frac{105}{8}-\frac{57}{2}\log 2 -\frac{33}{8}A_{0,2,0}-\frac{7}{4}A_{0,3,0}+\frac{3}{8}A_{0,4,0}; \\
A_{3,4,0} &=& -\frac{4211}{432}-\frac{769}{36}\log 2 -\frac{155}{48}A_{0,2,0}-\frac{37}{24}A_{0,3,0}+\frac{5}{16}A_{0,4,0}.
\end{eqnarray}
\end{exa}
The ASCII file in the \texttt{anc} directory contains a table with columns of $m$, $l$, $n$ and $A_{m,l,n}$.
A fifth column with $B_{m,l}$ is attached in the rows where $n=0$ and $l\ge 1$.

\section{Summary}
The index triple of $m,l,n\ge 0$ is covered
by the following algorithm to calculate $A_{m,l,n}$: 
\begin{itemize}
\item
Use \eqref{eq.recn2} and \eqref{eq.recn1} to reduce $n$ to zero;
\item
if $l=0$, use eventually \eqref{eq.l0};
\item
if $l=1$, use eventually \eqref{eq.l1};
\item
if $l\ge 2$, use \eqref{eq.recm} to reduce $m$ to zero and plug in the constants of Section \ref{sec.m0}.
\end{itemize}

\bibliographystyle{amsplain}
\bibliography{all}

\providecommand{\bysame}{\leavevmode\hbox to3em{\hrulefill}\thinspace}
\providecommand{\MR}{\relax\ifhmode\unskip\space\fi MR }
\providecommand{\MRhref}[2]{%
  \href{http://www.ams.org/mathscinet-getitem?mr=#1}{#2}
}
\providecommand{\href}[2]{#2}
\begin{thebibliography}{10}

\bibitem{AS}
Milton Abramowitz and Irene~A. Stegun (eds.), \emph{Handbook of mathematical
  functions}, 9th ed., Dover Publications, New York, 1972. \MR{0167642}

\bibitem{AdamchikJCAM79}
V.~S. Adamchik, \emph{On {S}tirling numbers and {E}uler sums}, J. Comput. Appl.
  Math \textbf{79} (1997), no.~1, 119--130. \MR{1437973}

\bibitem{BeumerAMM68}
M.~G. Beumer, \emph{Some special integrals}, Am. Math. Montly \textbf{68}
  (1961), no.~7, 645--647. \MR{0132150}

\bibitem{BowmanJLMS22}
F.~Bowman, \emph{Note on the integral $\int_0^{\pi/2}
  (log\sin\theta)^nd\theta$}, J. Lond. Math. Soc. (1947), no.~22, 172--173.
  \MR{0024493}

\bibitem{BoydAMC29}
John~P. Boyd, \emph{The asymptotic chebyshev coefficients for functions with
  logarithmic endpoint singularities: mappings and singular basis functions},
  App. Math. Comput. \textbf{29} (1989), no.~1, 49--67. \MR{0973493}

\bibitem{ChaudhuriAMM74}
Jyoti Chaudhuri, \emph{Some special integrals}, Am. Math. Monthly \textbf{74}
  (1967), no.~5, 545--548. \MR{0213613}

\bibitem{Clenshaw1954}
C.~W. Clenshaw, \emph{Polynomial approximations to elementary functions},
  Math.\ Tabl.\ Aids Comput. \textbf{8} (1954), no.~47, 143--147. \MR{0063487}

\bibitem{sloane}
O.~E. I.~S. Foundation~Inc., \emph{The {O}n-{L}ine {E}ncyclopedia {O}f
  {I}nteger {S}equences},  (2024), https://oeis.org/. \MR{3822822}

\bibitem{GottschalkJPA21}
J.~E. Gottschalk and E.~N. Maslen, \emph{Reduction formulae for generalised
  hypergeometric functions of one variable}, J. Phys. A: Math.\ Gen.
  \textbf{21} (1988), 1983--1998. \MR{0952917}

\bibitem{GR}
I.~Gradstein and I.~Ryshik, \emph{Summen-, {P}rodukt- und {I}ntegraltafeln},
  1st ed., Harri Deutsch, Thun, 1981. \MR{0671418}

\bibitem{KolbigMathComp40}
K.~S. K\"olbig, \emph{On the integral $\int_0^{\pi/2} \log^n\cos x\log^p \sin x
  dx$}, Math.\ Comp. \textbf{40} (1983), no.~162, 565--570. \MR{0689472}

\bibitem{KrupnikovJCAM78}
Ernst~D. Krupnikov and K.~S. K\"olbig, \emph{Some special cases of the
  generalized hypergeometric function $_{q+1}f_q$}, J. Comp. Appl. Math.
  \textbf{78} (1997), 79--95. \MR{1436781}

\bibitem{DLMF}
{Natl.\ Inst.\ Stand.\ Technol.}, \emph{Digital library of mathematical
  functions}, {NIST}, 2024. \MR{1990416}

\bibitem{QiRACSAM114}
Feng Qi and Chuan-Jun Huang, \emph{Computing sums in terms of beta, polygamma,
  and gauss hypergeometric functions}, Rev. Real Acad. Cien. Ex., Fis. Natur.
  Ser. A \textbf{114} (2020), 191.

\end{thebibliography}

\end{document}